\documentstyle{article}

\begin{document}

\newcommand{\iu}{\mbox{${{{\cal U}}}$}}
\newcommand{\iv}{\mbox{${{{\cal V}}}$}}
\newcommand{\ip}{\mbox{${{{\cal P}}}$}}
\newcommand{\ib}{\mbox{${{{\cal B}}}$}}
\newcommand{\xy}{(X,\iu,Y,/tau}
\newcommand{\qu}{quasi-uniformity}
\newcommand{\qus}{quasi-uniformities}
\newcommand{\qp}{quasi-proximity}
\newcommand{\vege}{\mbox{$\bullet$}}
\newcommand{\ce}{compatible extension}
\newcommand{\ut}{(\iu,\tau')}
\newcommand{\s}{\mbox{$\sigma$}}
\newcommand{\es}{\mbox{$\emptyset$}}
\renewcommand{\a}{{{\ \rm and\ }}}
\newcommand{\st}{{{\ \rm such\ that\ }}}
\renewcommand{\t}{{{\ \rm then\ }}}
\newcommand{\h}{{{\ \rm hence\ }}}
\renewcommand{\i}{{{\ \rm if\ }}}
\newcommand{\so}{{{\ \rm so\ }}}
\newcommand{\sa}{\ and\ }
\newcommand{\sst}{\ such\ that\ }
\newcommand{\sfe}{\ for\ every\ }
\newcommand{\sw}{\ where\ }
\newcommand{\si}{\ if\ }
\newcommand{\sthi}{\ there\ is\ }
\newcommand{\stha}{\ there\ are\ }
\newcommand{\nes}{\not=\es}
\newcommand{\w}{{{\ \rm where\ }}}
\renewcommand{\ss}{\subseteq}
\newcommand{\ki}{\backslash}
\newcommand{\fe}{{{\ \rm for\ every\ }}}
\newcommand{\fepymx}{\fe\pymx}
\newcommand{\sfepymx}{\sfe\pymx}
\newcommand{\il}{{\mbox{$\cal L$}}}
\newcommand{\im}{{{\ \rm implies\ }}}
\newcommand{\nss}{\not\ss}
\newcommand{\thi}{{{\ \rm there\ is\ }}}
\newcommand{\tha}{{{\ \rm there\ are\ }}}
\renewcommand{\P}{\par \rm Proof:\ }
\newcommand{\er}{\mbox{$\rm I\!\bf R$}}
\newcommand{\en}{\mbox{$\rm I\!\bf N$}}
\newcommand{\lett}{{{\rm \ let\ }}}
\newcommand{\ze}{{{\rm\bf Z\!\!\!\!\it Z}}}
\newcommand{\pni}{\par\noindent}
\newcommand{\filvx} {{\rm Fil}'(X)}
\newcommand{\ir}[1]{\mbox{${{{\cal #1}}}$}}
\newcommand{\iutp}{\mbox{$\iu^{tp}$}}
\newcommand{\filx}{\mbox{\rm fil$_X$}}
\newcommand{\filxx}{\mbox{\rm fil$_{X\times X}$}}
\newcommand{\expx}{\mbox{${\rm exp}X$}}
\newcommand{\pde}{\mbox{$\pi(\delta)$}}
\newcommand{\pdn}{\mbox{$\pi(\delta^0)$}}
\newcommand{\ivo}{\mbox{$\iv_{\omega}$}}
\newcommand{\ivd}{\mbox{$\iv_{\delta}$}}
\newcommand{\ivn}{\mbox{$\iv^0$}}
\newcommand{\ivs}{\mbox{$\iv_{\sigma}$}}
\newcommand{\Nt}{\mbox{$N(\tau)$}}
\newcommand{\Tt}{\mbox{$T(\tau)$}}
\newcommand{\NX}{\mbox{$N(X)$}}
\newcommand{\TX}{\mbox{$T(X)$}}
\newcommand{\ua}{\mbox{$U_{\alpha}$}}
\newcommand{\kaan}{\mbox{$2^{\aleph_0}$}}
\newcommand{\kakaan}{\mbox{$2^{2^{\aleph_0}}$}}
\newcommand{\Tk}{\mbox{\rm T${}_2$}}
\newcommand{\ipoc}{\mbox{interior preserving open cover}}

\newtheorem{th}{Theorem}[section]
\newtheorem{cor}[th]{Corollary}
\newtheorem{ex}[th]{Example}
\newtheorem{rem}[th]{Remark}
\newtheorem{prop}[th]{Proposition}
\newtheorem{lem}[th]{Lemma}
\newtheorem{defi}[th]{Definition}
\newtheorem{prob}{Problem}

\title{On the cardinality of \pde}
\author{A. Losonczi}
\date{March 26, 1999}
\maketitle

\begin{abstract}
\noindent
We prove that the cardinality of transitive \qus\ in a \qp\ class is at least
\kakaan\ if there exist at least two transitive \qus\ in the class. The
transitive elements of 
\pde\ are characterized if \ivd\ is transitive, and in this case we give a
condition when there exists a unique transitive \qu\ in \pde. 
\end{abstract}

\pni
AMS {\it Subject Classification:} 54E15; 54A25


\section{Introduction}

The purpose of this paper is the generalization of a result of \cite{c3} which
states that there exist either a unique or at least \kakaan\ compatible
transitive \qus\ on a topological space. Moreover we proved there that
$|\pi(\delta^1)|=1$ or $\ge\kakaan$ where $\delta^1$ denotes the finest
compatible \qp\ on $X$, and $\pi(\delta^1)$ is its \qp\ class, namely
$\pi(\delta^1)=\{\iv:\iv\supset\ip\}$ where \ip\ is the Pervin \qu\ of
$X$. The natural question arises : is something similar true for \pde\ if
$\delta$ is an arbitrary compatible \qp? In this paper we partly answer this
question by proving that $|\pde|=1$ or $\ge\kakaan$ if the coarsest element of
\pde\ is transitive.

In \cite{c8} we prove that in the class of infinite locally compact T${}_2$
spaces $|\pi(\delta^0)|=1$ if and only if $X$ is compact or non-Lindel\"of
($\delta^0$ denotes the coarsest compatible \qp) and if $X$ is non-compact and 
Lindel\"of, then $|\pi(\delta^0)|\ge\kakaan$. 

We give two elementary definitions.

\begin{defi} \label{nt} Let $(X,\tau)$ be a topological space. Then \NX\ or
$N(\tau)$ (\TX\ or $T(\tau)$)
denotes the set of all compatible (transitive) \qus\ on $X$ respectively. 
\end{defi}

\begin{defi} A base \ib\ of a topological space is called an l-base or a
lattice-base if it is
closed under finite union and finite intersection and\
$\es,X\in\ib$. \end{defi} 

We enumerate some results connected with this notion. For the proofs of these
results the reader may wish to consult \cite{c2}.

We say that $\alpha$ is an interior preserving open cover if it is an open
cover and for every $x\in X$ the set $\bigcap\{N\in\alpha:x\in N\}$ is
open or equivalently $\alpha'\ss\alpha\im\bigcap\alpha'$ is open. We can
assign a transitive neighbournet \ua\ to $\alpha$ in the 
following way: $\ua(x)=\bigcap\{N\in\alpha:x\in N\}\ (x\in X)$. If
$\alpha'=\{\ua(x):x\in X\}\t\alpha'$ is also an \ipoc\a$U_{\alpha'}=\ua$. An
important remark is that if $\alpha$ is finite then $\{\ua(x):x\in X\}$ is
finite too and if $\ua\in\ip\t\alpha$ is finite where \ip\ denotes the Pervin
\qu\ of $X$ namely $\ip=\filxx\{\ua:\alpha$ is finite$\}$. It is known that a
\qu\ is totally bounded if and only if it is contained in \ip\ (see
\cite{fl}1). 

In \cite{c2} we showed that there is a one-to-one correspondence between the
set of compatible totally bounded transitive \qus\a the set of l-bases. Namely
if \iv\ is a \qu\ of the mentioned type then $\ib(\iv)=\{N\in\tau:U_N\in\iv\}$
is an l-base, and if \ib\ is an l-base then $\iv(\ib)=\filxx\{U_N:N\in\ib\}$
is a totally bounded transitive \qu, where $U_N=(N\times N)\cup ((X-N)\times
X)=U_{\{N,X\}}$ for the cover$\{N,X\}$. 

If we say that $\delta$ is a
quasi-proximity we use it in the sense of \cite{fl}1.22. If \iv\ is a \qu,
then $\delta(\iv)\a\tau(\iv)$ will always denote the quasi-proximity and the
topology induced by \iv\ respectively, $\tau(\delta)$ will denote the
topology induced by $\delta$. Let $\pde=\{\iv:\delta(\iv)=\delta\}$. We know
from \cite{fl}1.33 that \fe $\delta,\ \pde\nes$. Moreover there exists a
coarsest element of \pde, it is denoted by $\iv_{\delta}$. This is totally
bounded and the only totally bounded member of \pde. If a \qu\ \iv\ is given
then $\ivo$ denotes the coarsest element of $\pi\bigl(\delta(\iv)\bigr)$.


\section{The main results}

First we want to characterize the transitive elements of \pde\ where $\delta$
is a compatible \qp\st \ivd\ is transitive. To this aim we will need some
lemmas. 

\begin{lem}\label{l01} Let $\delta$ be a quasi-proximity on $X$ compatible
with $\tau$ 
such that $\iv_{\delta}$ is transitive. Let $\iv\in T(X)\sst \ivd\ss\iv$. In
this case $\iv\in\pde$ if and only if $N\in\tau,\ U_N\in\iv$ imply that
$N\in\ib(\delta)=\ib(\ivd)$. \end{lem} 

\P By \cite{c1}2.6 we know that $(\ivo)_t=(\ip)_t\cap (\iv)_t$ where for a \qu\
\iu,\ $(\iu)_t=\{U\in\iu:U$ is transitive$\}$. 

If $\iv\in\pde$ then
$\ivo=\ivd$. By the previous observation if $U_N\in\iv\t U_N\in\ivd\h
N\in\ib(\delta)$. Let us prove the sufficiency. We get
$(\iv)_t\cap(\ip)_t\ss\ivd\a \ivo\ss\ivd$ which yields that
$\iv\in\pde$. \vege     

\begin{lem} \label{l1} Let $\iv\in\Tt,\ \ib=\ib(\ivo)\sa U\in\iv$ be
transitive. In this case $U(A)\in\ib$ if $A\ss X$. \end{lem}

\P Let $A\ss X\a N=U(A)$. It is easy to check that $U\ss U_N$ so $U_N\in\iv,\a
U_N\in\ivo$.
\vege

\begin{cor} \label{cor1} Let $\iv\in\Tt,\ \ib=\ib(\ivo)\sa \ua\in\iv$ where
$\alpha$ is an interior preserving open cover. Then $A\in\alpha$ implies that
$A\in\ib$. \end{cor} 

\P $A=\ua(A)$. \vege

\begin{prop} \label{p1} Let $\delta$ be a compatible quasi-proximity on $X$
such that $\iv_{\delta}$ is transitive and let \ib\ be the l-base associated
with \ivd.   
Let $\iv=\filxx\{\ivd,U_i:i\in I\}$ where $U_i$ is transitive ($i\in I$) and
the system $\{U_i:i\in I\}$ is closed under finite intersection. Then
$\iv\in\pde$ if and only if $\forall i\in I,\ \forall A\ss X\ U_i(A)\in\ib$. 
\end{prop}

\P If $\iv\in\pde\t\ivo=\ivd$ and we get the statement by the previous lemma
(\ref{l1}). 

To prove the opposite case let $N\in\tau=\tau(\delta),\ U_N\in\iv$. We need
that $N\in\ib$ by \ref{l01}. We know 
that there are $M_i\in\tau,\ U_{M_i}\in\ivd\a j\in I\st$ 
\[\bigcap_{i=1}^n U_{M_i}\cap U_j\ss U_N\]
or in other words $U_{\{M_1,\ldots,M_n,X\}}\cap U_j\ss U_N$ for the cover
$\{M_1,\ldots,M_N,X\}$. Obviously 
$M_i\in\ib$ and we can assume the system $\{M_1,\ldots,M_n,X\}$ is closed
for union and intersection and $M_n=X$. Let $U'=U_{\{M_1,\ldots,M_n\}}$.
If $x\in N\t\thi M_i\st x\in M_i\cap U_j(x)\ss U_N(x)=N$. Now
\[N=\bigcup_{x\in N} (U_{\{M_1,\ldots,M_n\}}\cap U_j)(x)=\]

\[\bigcup_{x:U'(x)=M_1}(M_1\cap U_j(x))\bigcup\bigcup_{x:U'(x)=M_2}(M_2\cap
U_j(x))\bigcup\ldots\bigcup \bigcup_{x:U'(x)=M_n}(M_n\cap
U_j(x))=\]

\[\bigl(M_1\cap U_j(\{x:U'(x)=M_1\})\bigr)\bigcup
\ldots\bigcup\bigl(M_n\cap U_j(\{x:U'(x)=M_n\})\bigr)\in\ib\]
since $U_j(\{x\in N:U'(x)=M_k\})\in\ib$ by assumption and \ib\ is 
an l-base. \vege

\begin{cor} \label{fc} Let $\delta$ be a quasi-proximity, \ivd\ be transitive,
\ib\ be the l-base associated with it and \iv\ be a transitive \qu\ on $X$.Then
$\iv\in\pde$ if and only if $\ivd\ss\iv\sa$for every transitive $U\in\iv$ and
for every $A\ss X,\ U(A)\in\ib$. \end{cor}

\P This is an obvious consequence of \ref{p1}. \vege
 
Now we can give condition for $|\pde\cap\TX|>1$.


\begin{th} \label{t1} Let $\delta$ be a compatible quasi-proximity,\ \ivd\ 
be transitive and \ib\ be
the l-base associated with \ivd. In this case there exists 
$\iv\in\pde\cap\Tt,\ \iv\not=\ivd$ if and only if either

1. there exists a system of sets $\{N_i:i\in \en\}\sst N_i\in\ib,
N_i\ss N_{i+1}, N_i\not=N_{i+1}\sa \bigcup_1^{\infty}N_i\in\ib$ or

2. there exists a system of sets $\{N_i:i\in \en\}\sst N_i\in\ib,
N_{i+1}\ss N_i, N_i\not=N_{i+1}\sa \bigcap_1^{\infty}N_i\in\ib$.
\end{th}

\P Let us prove first the sufficiency and suppose that condition 1 holds. Let
$\alpha=\{X;N_i:i\in\en\}\a\iv=\filxx\{\ivd,U_{\alpha}\}$. By \ref{p1}
$\iv\in\pde\cap\Tt$ since if $A\ss X$ then $U_{\alpha}(A)$ is either a finite
union of $N_i$'s or if it is infinite then it equals to
$\bigcup_1^{\infty}N_i\in\ib$. If 
2 holds then let $\alpha=\{X;\bigcap_1^{\infty}N_i\}\cup\{N_i:i\in\en\}\a\iv=
\filxx\{\ivd,U_{\alpha}\}$. Then $U_{\alpha}(A)$ is always a finite union so
it is in \ib. Finally by \cite{c2}2.7 $\ua\notin\ivd$ in both cases.

Let us prove the necessity. If $\iv\in\pde$ is transitive, $\iv\not=\ivd$ then
there 
is a transitive $U\in\iv\st\{U(x):x\in X\}$ is infinite. If we
suppose that there is $A=\{x_i\in X:i\in\en\}\st$ 
$U(x_n)\nss\bigcup_{i=1}^{n-1}U(x_i)\ (\forall n\in\en)$ then condition 1
holds since by \ref{fc} $U(x_i)\in\ib,\  
U(A)\in\ib$ and the system $\{\cup_{j=1}^nU(x_j):n\in\en\}$ will work. 

In case there exists no such system then there is a finite set $A_1\ss X$
such that $X=U(A_1)$. Since $\{U(x):x\in X\}$ is infinite then there exists
$x_1\in A_1$ such that $\{U(y):y\in U(x_1)\}$ is infinite by $U$ being
transitive. Then there is a finite $A_2\ss X$ such that $U(A_2)=U(x_1)$ and
there exists $x_2\in A_2\st \{U(y):y\in 
U(x_2)\}$ is infinite. Let us continue this process. We get a system
$\{U(x_i):i\in\en\}\st U(x_i)\supset U(x_{i+1})$. $U(x_i)\not=U(x_{i+1})$ can
also be assumed. 
By \ref{cor1} $\bigcap_1^{\infty} U(x_i)\in\ib$. 
\vege

\begin{cor} \label{coa} If $\delta_1$ is coarser than $\delta_2$,
$\iv_{\delta_1},\iv_{\delta_2}$ are transitive, $|\pi(\delta_1)\cap\Tt|>1$
then $|\pi(\delta_2)\cap\Tt|>1$. \end{cor}

\P Obviously $\ib(\delta_1)\ss\ib(\delta_2)$ and apply \ref{t1}. \vege

%
%
%

We will need definitions and a theorem from \cite{c3} dealing with p-filters.

\begin{defi} (\cite{c3}2.2) Let $N\ss\en$. We call a subset $H$ of $N$ a pile
of $N$ or an 
$N$-pile or simply pile (if there is no misunderstanding) if
\par 1. $i,j\in H\sa i<k<j,\ k\in\en$ implies that $k\in H$ and
\par 2. $H$ is maximal for the previous property.\end{defi}

\begin{defi} 
(\cite{c3}2.3) 
Let again $N\ss\en$. We say that $Z\ss\en$ is
admissible for $N$ 
if there exists $k\in\en\sst H$ being an $N$-pile implies that $|H\cap Z|\le
k$. \end{defi}

\begin{defi} (\cite{c3}2.4) We call a filter \s\ on \en\ a p-filter if,
whenever $N\in\s\sa 
Z$ is admissible for $N$, then $N-Z\in\s$. \end{defi}

\begin{th} (\cite{c3}2.7) $|\{\s:\s$ is a p-filter
$\}|=2^{2^{\aleph_0}}$. \end{th} 

Now we are ready to prove the main theorem.

\begin{th} \label{fotetel} Let $\delta$ be a \qp\ such that \ivd\ is
transitive. In this case $|\pde|=1$ or $|\pde|\ge 2^{2^{\aleph_0}}$, moreover
$|\pde\cap\Tt|=1$ or $\ge 2^{2^{\aleph_0}}$. \end{th}

\P Let \ib\ be the l-base associated with \ivd.  
If $|\pde|>1$ then condition 1 or 2 holds in theorem \ref{t1}.

1. Suppose that there exists a system of sets $\{N_i:i\in \en\}\st N_i\in\ib, N_i\ss N_{i+1}, N_i\not=N_{i+1}\a \bigcup_1^{\infty}N_i\in\ib$.

If $A\ss\en\t$ let
\[\alpha_A=\{X,N_i:i\notin A\}.\] Then $\alpha_A$ is an interior preserving  
open cover of $X$. Let $U_A(x)=U_{\alpha_A}(x)=\cap\{M\in\alpha_A:x\in M\}$.
If $\sigma$ is a p-filter on $\en$ then let
\[\iv_{\sigma}=\filxx\{\ivd,U_A:A\in\sigma\}.\] 
It is obvious that
$\iv_{\sigma}\in\Tt$. It is also straightforward that $U_{A_1}\cap
U_{A_2}=U_{A_1\cap A_2}$. By \ref{p1}, $\ivs\in\pde$. 

  
We show that if
$\sigma_1\not=\sigma_2\t\iv_{\sigma_1}\not=\iv_{\sigma_2}$. This yields
immediately that the cardinality of the set of all p-filters on $\en$ is less
than or equal to the cardinality of $\pde\cap\Tt$, hence $|\pde\cap\Tt|\ge
2^{2^{\aleph_0}}$. 

Now let $\iv_{\sigma_1}\ss\iv_{\sigma_2}$ then we will show that
$\sigma_1\ss\sigma_2$, which is enough to prove since with the opposite 
case it implies the required statement.

Let $A_1\in\sigma_1$. Then $U_{A_1}\in\iv_{\sigma_1}$ hence \tha
$A_2\in\sigma_2\a P=U_{\beta}\in\ivd\st U_{\beta}\cap U_{A_2}\ss U_{A_1}$ where
$\beta$ is a finite subset of \ib. 
Let $|\{P(x):x\in X\}|=k\in\en$. Let $H\ss A_2$ be an
$A_2$-pile\st$H=\{r\in\en:p<r<q\}$. Suppose that $i\in H-A_1$ and let us fix 
$x_i\in N_i-N_{i-1}$. Then $U_{A_1}(x_i)=N_i\a U_{A_2}(x_i)=N_q$ and we get
that $N_q\cap 
P(x_i)\ss N_i$ which implies that $P(x_i)\cap(N_q-N_i)=\es$. If $j\in H-A_1\st
i<j\t x_j\in P(x_j)\cap (N_q-N_i)\nes$ so $P(x_i)\not=P(x_j)$ which verifies
that the function $f$ defined by $f(i)=P(x_i)\ (i\in H-A_1)$ is
injective. A similar argument applies if $q=\infty$. Hence $|H-A_1|\le 
k$ and there is $Z=A_2-A_1\ss\en\st A_2-Z\ss A_1$ where $Z$ is admissible for
$A_2$. Since $\s_2$ is a p-filter then $A_1\in\s_2\a\s_1\ss\s_2$.

2. Suppose that there exists a system of sets $\{N_i:i\in \en\}\st N_i\in\ib,
N_{i+1}\ss N_i, N_i\not=N_{i+1}\a \bigcap_1^{\infty}N_i\in\ib$.

Let $N_0=X$. If $A\ss\en\t$ let
\[\alpha_A=\{X,\bigcap_{j=1}^{\infty}N_j\}\bigcup\{N_i:i\notin A\}.\] 
Then $\alpha_A$ is an interior preserving open cover of $X$. Let
$U_A(x)=\bigcap\{M\in\alpha_A:x\in M\}$. If $\sigma$ is a p-filter on $\en$
then let 
\[\iv_{\sigma}=\filxx\{\ivd,U_A:A\in\sigma\}.\]
It is obvious that
$\iv_{\sigma}\in\Tt$. It is also straightforward that $U_{A_1}\cap
U_{A_2}=U_{A_1\cap A_2}$. By \ref{p1}, $\ivs\in\pde$. 

We show that if
$\sigma_1\not=\sigma_2\t\iv_{\sigma_1}\not=\iv_{\sigma_2}$. This yields
immediately that $|\pde\cap\Tt|\ge 2^{2^{\aleph_0}}$. 

Now let $\iv_{\sigma_1}\ss\iv_{\sigma_2}$ then we will show that
$\sigma_1\ss\sigma_2$, which is enough to be proved.

Let $A_1\in\sigma_1$. Then $U_{A_1}\in\iv_{\sigma_1}$ hence \tha
$A_2\in\sigma_2\a P=U_{\beta}\in\ivd\st U_{\beta}\cap U_{A_2}\ss U_{A_1}$ where
$\beta$ is a finite subset of \ib. 
Let $|\{P(x):x\in X\}|=k\in\en$. Let $H\ss A_2$ be an
$A_2$-pile\st$H=\{r\in\en:p<r<q\}$. (Possibly $p=0$ or $q=\infty$.) Suppose
that $i\in H-A_1$ and let us fix  
$x_i\in N_i-N_{i+1}$. Then $U_{A_1}(x_i)=N_i\a U_{A_2}(x_i)=N_p$ and we get
that $N_p\cap  
P(x_i)\ss N_i$ which implies that $P(x_i)\cap(N_p-N_i)=\es$. If $j\in H-A_1\st
i<j\t x_j\in P(x_j)\cap (N_p-N_i)\nes$ so $P(x_i)\not=P(x_j)$ which verifies
that the function $f$ defined by $f(i)=P(x_i)\ (i\in H-A_1)$ is
injective. Hence $|H-A_1|\le 
k$ and there is $Z\ss\en\st A_2-Z\ss A_1$ where $Z$ is admissible for
$A_2$. Since $\s_2$ is a p-filter then $A_1\in\s_2\a\s_1\ss\s_2$. \vege  

We know from \cite{c8}2.4 that in the class of locally compact, \Tk\ spaces 
if there is a
transitive $\iv\in\pdn$ then $X$ is 0-dimensional. With the help of \ref{t1}
one can easily prove that $|\pdn\cap\Tt|=1$ if $X$ is compact or
non-Lindel\"of,  but \cite{c8}2.14 and 3.12 yield a more general result,
namely $|\pdn|=1$ in these cases. 
So in the context of this paper only the case in which 
$X$ is non-compact, Lindel\"of and 0-dimensional
is interesting.

\begin{cor} Let $X$ be locally compact, \Tk, non-compact, Lindel\"of and
0-dimensional. If $\delta$ is a compatible \qp\sst $\ivd$ is transitive then
$|\pi(\delta)\cap\Tt|\ge\kakaan$. \end{cor}

\P By \ref{coa} it is enough to verify the statement for $\delta^0$. In this
case $\ib=\{$compact-open sets$\}\cup\{\es,X\}$, and there is a strictly
increasing sequence $(N_i)$ of compact-open sets,\st $\cup_1^{\infty}N_i=X$ and
\ref{t1} and \ref{fotetel} are applicable. \vege

Finally we present some open problems which seem to be interesting. 

\begin{prob} What can we say about $|\pde|$ if \ivd\ is not transitive?
\end{prob} 

Remark : We note that in the meantime this problem
have been solved by K\"unzi in \cite{kL}.

\begin{prob} Let \ivd\ be transitive. 
What can be said about $|\pde\cap(N(X)-T(X))|$?
Can it occur that there exists no non-transitive \qu\ in \pde\ if
$|\pi(\delta)|>1$? \end{prob}

\begin{prob}Is there any connection between
$|\pi(\delta_1)|\sa|\pi(\delta_2)|$ if $\delta_1$ is coarser than $\delta_2$?
Is it true that $|\pi(\delta_1)|\le|\pi(\delta_2)|$ in this case? \end{prob}

\begin{prob}Is there always a non-transitive \qu\ that is finer than \ip
(assuming that $|\pi(\delta^1)|>1$)? (see also \cite{kuj}Remark 1)
\end{prob} 



\smallskip

\noindent
Alfr\'ed R\'enyi Institute of Mathematics, Hungarian Academy of Sciences,

\noindent Re\'altanoda u. 13-15, H-1364 Budapest, Hungary 
\smallskip

\noindent
E-mail: losonczi@math-inst.hu

\end{document}